\documentclass[a4paper,10pt]{article} \usepackage{amsmath,amssymb,amsthm}
\usepackage[english]{babel}

\newtheorem{theo}{Theorem}[section] 
\newtheorem{defi}[theo]{Definition}
\newtheorem{lemm}[theo]{Lemma} 
\newtheorem{prop}[theo]{Proposition}
\newtheorem{coro}[theo]{Corollary}

\newcommand{\Na}{\mathbb N}                   

\newcommand{\Ra}{\mathbb R}                   

\newcommand{\Ca}{\mathbb C}                   

\newcommand{\scal}[1]{\langle #1 \rangle}

\newcommand{\finpreuve}{\hfill $\Box$}

\newcommand{\name}{$\underline{\qquad \qquad}$}

\begin{document}

\author{  Jean-Marc Bouclet}
\title{Low frequency estimates for long range perturbations in divergence form}

\maketitle

\begin{abstract} We prove a uniform control as $ z \rightarrow 0 $ for the  resolvent $ (P-z)^{-1} $  of long range perturbations $  P $ of the Euclidean Laplacian  in divergence form.  
\end{abstract}

\section{Introduction and main results}
Consider an elliptic self-adjoint operator in divergence form on $ L^2 (\Ra^d) $, $ d \geq 2 $,
$$ P = - \mbox{div} \left( G (x) \nabla \right) , $$
where $ G(x) $ is a $ d \times d $ matrix with real entries satisfying, for some $ \Lambda_+ \geq \Lambda_- > 0 $, 
\begin{eqnarray}
 G (x)^T = G (x), \qquad \Lambda_+ \geq G (x) \geq \Lambda_-  , \qquad x \in \Ra^d . \label{bornesellipticite}
\end{eqnarray}
 Throughout the paper, we shall assume that $ G $ belongs to $ C^{\infty}_b (\Ra^d) $ ie that $ \partial^{\alpha} G $ has bounded entries for all multiindices $ \alpha $, but this is mostly for convenience and much weaker assumptions on the regularity of $ G $ could actually be considered. For instance, in polar coordinates $ x = |x|\omega $,  Theorem \ref{theoremeprincipal} below will not use any regularity in the angular variable $ \omega $.

We mainly have in mind  long range perturbations of the Euclidean Laplacian, namely the situation where, for some $ \mu > 0 $,
\begin{eqnarray}
\big| \partial^{\alpha} \left( G (x) - I_d \right) \big| \leq C_{\alpha} \scal{x}^{-\mu-|\alpha|}, \qquad x \in \Ra^d, \label{definitlongueportee}
\end{eqnarray}
$ I_d $ being the identity matrix and $ \scal{x} = (1+|x|^2)^{1/2} $ the usual japanese bracket. In this case, it is well known that the resolvent
$ ( P - z )^{-1} $ satisfies the limiting absorption principle, ie that the limits
$$ ( P - \lambda \mp i  0 )^{-1} := \lim_{\delta \rightarrow 0^+} (P - \lambda \mp i \delta)^{-1} $$
exist at all positive energies $ \lambda > 0 $ (the frequencies being $ \lambda^{1/2} $) in weighted $ L^2 $ spaces (see the historical papers \cite{Agmon,Mourre}, the references therein and the references below on quantitative bounds).
 Typically, for all $ \lambda_2 > \lambda_1 > 0 $ and all $ s > 1/2 $, we have bounds of the form
\begin{eqnarray}
 \big| \big| \scal{x}^{-s} (P- \lambda - i 0)^{-1} \scal{x}^{-s} \big| \big|_{L^2 \rightarrow L^2} \leq C (s,\lambda_1,\lambda_2), \qquad \lambda \in [\lambda_1,\lambda_2 ] , \label{absorptionprinciple}
\end{eqnarray}
and the same holds of course for $ (P- \lambda + i 0)^{-1} $ by taking the adjoint. The behaviour of the constant $ C (s,\lambda_1,\lambda_2) $ is very well known as long as $ \lambda_1 $ doesn't go to $ 0 $. For a fixed energy window, the results follow essentially from the Mourre theory \cite{Mourre} since one knows that there are no embedded eigenvalues for such operators \cite{KochTataru}. At large energies, $ \lambda_1 \sim \lambda_2 \rightarrow \infty $, $ C (s , \lambda_1 , \lambda_2) $ is at worst of order $ e^{C \lambda_2^{1/2}} $, see \cite{CardosoVodev}, but can be taken of order $ \lambda_1^{-1/2} $ if there are no trapped geodesics  (ie all geodesics escape to infinity) - see\cite{GerardMartinez,Wang0,Robert,Burq2,VasyZworski}. 

\medskip

In this paper, we address the problem of the behaviour of such estimates as $ \lambda_1  \downarrow 0$. Let us recall that a quick look at the kernel of the resolvent in the flat case ($ P = - \Delta $), whose kernel is given for $ d = 3 $ (for simplicity) by
\begin{eqnarray}
 K_{\rm flat}(x,y,z) = \frac{e^{i z^{1/2}|x-y|}}{4\pi|x-y|}, \qquad \mbox{Im}(z^{1/2}) \geq 0 , \label{noyauLaplacien}
\end{eqnarray}
suggests that, if one has no oscillation, ie if $ z = 0 $, one should rather choose $ s > 2 $ than $ s > 1/2 $ in (\ref{absorptionprinciple}). This (natural) restriction is however essentially irrelevant for us: our point in the present paper is not to get the sharpest weights (e.g. work in optimal Besov spaces) but only to get a control on $ \scal{x}^{-s}(P-\lambda-i0)^{-1} \scal{x}^{-s} $  as $ \lambda \rightarrow 0  $, for some  $ s $.

The very natural question of low frequency asymptotics for the resolvent of Schr\"odinger type operators has been considered in many papers. However the situation is not as clear as for the positive energies. For perturbations of the flat Laplacian by potentials, we refer to
 \cite{JensenKato,Yafaev,Murata,Nakamura,JensenNenciu,Wang1,FournaisSkibsted}, to the references therein  and also to the recent very detailled study \cite{DerezinskiSkibsted}.
 In a sense, perturbations by potentials are harder to study due to the possible resonances or (accumulation of) eigenvalues at $ 0 $. 
 
 For compactly supported perturbations of the flat Laplacian by  metrics and obstacles,  the behaviour of the resolvent at $ 0 $ is obtained fairly shortly in \cite{Morawetz,Burq1} but using strongly  the compact support assumption. Recently, Guillarmou and Hassell have investigated  carefully the low energy asymptotics of Schr\"odinger operators on asymptotically conical manifolds \cite{GuillarmouHassell,GuillarmouHassell2}. Using the sophisticated pseudo-differential calculus of Melrose, they are able to describe accurately the kernel of the Green function at low energies. In particular, they  derive optimal $ L^p $ bounds for the Riesz transform. This technology is also used in \cite{CarronCouhlonHassell}, again for the study of the range of $ p $ for which the Riez transform is $ L^p $ bounded. In a close geometric context, for very short range perturbations of exact conical metrics, Wang \cite{Wang2} also proves asymptotic expansion of the resolvent at low energies. 
 
 All the above papers dealing with metrics use a relatively strong decay of the perturbation at infinity  or assume at least certain asymptotic expansions which, in any case, exclude most long range perturbations. 
 
 The first message of the present paper is that nothing nasty can happen for long range perturbations of the metric. More precisely, we will show that, if the perturbation is uniformly small on $ \Ra^d $ (but arbitrarily long range at infinity), we have uniform bounds on the resolvent at low frequency. The second message is that, for arbitrary long range perturbations,  we can use certain properties of the Riesz transform to handle the non small compact part of the perturbation and get low energy estimates. In a sense, this is the opposite point of view to \cite{CarronCouhlonHassell,GuillarmouHassell,GuillarmouHassell2}, to the extent that we use the Riesz transform to analyze the resolvent instead of using information on the resolvent to study the Riesz transform. We furthermore think that the method described in this paper is quite simple (at least on $ \Ra^d $). More importantly, we hope that it is rather flexible and can be adapted to other  geometries.

Our main result is the following one.

\begin{theo} \label{theoremeconditionnel} Let $ d \geq 3 $. Assume that $ G  $ satisfies (\ref{bornesellipticite}) and (\ref{definitlongueportee}). Then, for all $ \epsilon > 0 $,
\begin{eqnarray}
 || \scal{x}^{-2-\epsilon} (P-z)^{-1} \scal{x}^{-2-\epsilon} ||_{L^2 \rightarrow L^2} \leq C_{\epsilon,G} , \qquad |z| < 1, \ z \notin \Ra . 
 \label{bassefrequencelongueportee}
\end{eqnarray}  
\end{theo}
In a sense, the proof of this theorem will be reduced to the case of a uniformly small perturbation of the Euclidean metric by "absorbing" the contribution of the (arbitrarily large) compact part of the perturbation using the local compactness of the Riesz transform $ \nabla P^{-1/2} $ (ie that $ \chi (x) \nabla P^{-1/2} \varphi (P) $ is compact for all $ \chi \in C_0^{\infty}(\Ra^d) $ and $ \varphi \in C_0^{\infty}(\Ra) $ - see Subsection \ref{soussectionRiesz} below).

\bigskip

The other results of this paper deal with small perturbations. To state them, we introduce the space $ S_{\rm dil} (\Ra^d) $ defined by
\begin{eqnarray}
 a \in S_{\rm dil}(\Ra^d) \qquad \Leftrightarrow \qquad a \in C^{\infty}_b (\Ra^d) \qquad \mbox{and} \qquad ( x \cdot \nabla )^n a \in L^{\infty} (\Ra^d) \ \ \mbox{for all }n , 
\end{eqnarray}
 and the related (semi-)norms
\begin{eqnarray}
  || a ||_{N,{\rm dil}} := \max_{n \leq N} || ( x \cdot \nabla )^n a ||_{L^{\infty}} . \label{seminormes}
\end{eqnarray}
For matrices $ H = (b_{jk}) $ with entries in $ S_{\rm dil}(\Ra^d) $, we shall denote $ || H ||_{N,{\rm dil}} $ for $ \max_{1 \leq j,k \leq d} || b_{jk} ||_{N,{\rm dil}} $.

As mentionned above, the condition  $ a \in C^{\infty}_b (\Ra^d) $ is mainly  for convenience, to simplify certain algebraic manipulations. For instance, it ensures that the resolvent $ (P-z)^{-1} $ maps the Schwartz space $ {\mathcal S}(\Ra^d) $ into itself if $ z \notin \Ra $, which is useful to compute commutators.

This space is obviously closely related to the following well known  generator of dilations,
\begin{eqnarray}
A = \frac{x \cdot \nabla + \nabla \cdot x}{2 i} = \frac{x \cdot \nabla}{i} + \frac{d}{2i} , \label{generateur}
\end{eqnarray}
so called for it is the self-adjoint generator of the unitary group on $ L^2 (\Ra^d) $ given by
\begin{eqnarray}
\left( e^{i t A} \varphi \right) (x) = e^{\frac{td}{2}} \varphi (e^t x) . \label{groupe}
\end{eqnarray}





\begin{theo} \label{theoremeprincipal} Assume that $ d \geq 2 $. Let $ G$ be of the form
\begin{eqnarray}
 G (x) = I_d + H (x) , \label{metrique}
\end{eqnarray}
with $ H $ symmetric and with real entries in $ S_{\rm dil} (\Ra^d) $. Then, for all $ \varepsilon >0 $, there exists $ C_{\varepsilon}  > 0 $  such that
for all $ H $ satisfying 
\begin{eqnarray}
2 G (x) - (x \cdot \nabla) H (x) \geq \varepsilon 
, \qquad x \in \Ra^d , \label{petitesse}
\end{eqnarray}
 and all $ h $ such that  $ 0<h \leq C_{\varepsilon}^{-1} (1+||H||_{4,{\rm dil}})^{-1}  $, we have
\begin{eqnarray}
 \left| \left| |D| (hA+i)^{-1} ( P - z )^{-1} (hA-i)^{-1} |D| \right| \right|_{L^{2} \rightarrow L^{2}} \leq \frac{C_{\varepsilon}}{h} 
 ,\qquad z \in \Ca \setminus \Ra . \label{inegaliteprincipale}
\end{eqnarray}
Here $ |D| $ is the usual Fourier multiplier by $ |\xi| $.
\end{theo}

The main point in this theorem is the uniform control in $z$ of the resolvent under the condition (\ref{petitesse}) which is essentially a smallness condition for it clearly holds if $ || H||_{1,{\rm dil}} $ is small enough. 
The main novelty is that we get bounds for small $z$, say $ |z| < 1 $. We however also obtain bounds for large $ z $ but these are essentially well known since the condition (\ref{petitesse}) implies that the metric $ G $ (or rather $ G^{-1} $) is non trapping ($ x \cdot \xi $ is a global escape function - see for instance \cite{GerardMartinez,Robert}). 

 We also give an explicit range of $h$ to show that the right hand side of (\ref{inegaliteprincipale}) can be estimated by $ (1+||H||_{4,{\rm dil}}) $ (for a suitable $h$) which could be of interest for instance for certain nonlinear applications. We point out that the regularity $ || H ||_{4,{\rm dil}} $ is probably not sharp. We have not tried to get the optimal regularity in order to avoid technicalities in the proofs and to focus on the main simple algebraic ideas; we thus might have done some relatively crude estimates at certain steps (in particular in Proposition \ref{integrationresolvente}). One may however hope to improve the regularity condition by changing $ ||H||_{4,{\rm dil}} $ into $ ||H||_{2,{\rm dil}} $.


\bigskip

We now derive weighted estimates of the same form as (\ref{absorptionprinciple}). 
For $ d \geq 3 $, recall the standard notation for the usual conjugate Sobolev exponents
$$ 2_* = \frac{2d}{d+2}, \qquad 2^* = \frac{2d}{d-2} . $$
\begin{coro} \label{corollaireTao} If $ d \geq 3 $, under the same assumptions as in Theorem \ref{theoremeprincipal}, we have
$$ \left| \left|  (hA+i)^{-1} ( P - z )^{-1} (hA+i)^{-1} \right| \right|_{L^{2_*} \rightarrow L^{2^*}} \leq \frac{C}{h}  ,\qquad z \in \Ca \setminus \Ra . $$
\end{coro}

This in turn leads to weighted estimates for long range perturbations of the Euclidean metric.
\begin{coro} \label{longueportee} Let $ d \geq 3 $. If $ G = I_d + H $ satisfies (\ref{bornesellipticite}), (\ref{definitlongueportee}) and (\ref{petitesse}), then
 for all $ \epsilon > 0 $, (\ref{bassefrequencelongueportee}) holds for all $ z \in \Ca \setminus \Ra $.
\end{coro}
Note that the difference between this corollary and Theorem \ref{theoremeconditionnel} is that the estimates hold for $z \in \Ca \setminus \Ra$. The latter is natural since the assumption (\ref{petitesse}) implies the non trapping condition, which gives the uniform control at high energies.

It is also worth noticing that the assumptions of Theorem \ref{theoremeprincipal} and the scale invariant space $ S_{\rm dil} (\Ra^d) $ are very close to the context of \cite{stTataglobal} where the time dependent Schr\"odinger equation is studied. Among other dispersive estimates,  Tataru proves in \cite{stTataglobal} $ L^2 $-space-time bounds, usually refered to as  global smoothing effect, for small long range perturbations of the euclidean metric, possibly time dependent, by using also positive commutator techniques. In the time independent case,  our (weighted) resolvent estimates (\ref{bassefrequencelongueportee}) combined with the usual ones at high energy also imply this smoothing effect. From the point of view of space-time bounds, the results of \cite{stTataglobal}  are stronger since they allow time dependent metrics. But on the other hand, in the time independent case, our resolvent estimates (which are $ L^{\infty}_{\rm loc} $ in term of the spectral parameter $z$) are stronger than $ L^2 $-space-time bounds on the evolution group. 

\bigskip




Rather vaguely, our results can be summarized as follows. Theorem \ref{theoremeprincipal} or Corollary \ref{longueportee} describe the situation near infinity, ie where one is close to the Laplacian. The intuition is that the Green function of the resolvent oscillates at large distances (consider typically  (\ref{noyauLaplacien})  for $ |x-y| \gtrsim \mbox{Re}(z^{-1/2}) $) which allows to use  commutator techniques. However, by the uncertainty principle, the energy localization close to $ 0 $ (ie classically at the point $ \xi =0 $) which is stronger than a localization at a positive energy, say close to $ 1 $ (ie on the sphere $ |\xi| = 1$), requires a stronger spatial delocalization which is responsible for the extra weights $ |D| $ in Theorem \ref{theoremeprincipal} or larger powers of $ \scal{x}^{-1} $ in Corollary \ref{longueportee}. On the other hand, at (relatively) short distances, one can not use oscillations. One is rather in a "singular integral" regime and, to this extent, the use of the Riesz transform is natural in the proof of Theorem \ref{theoremeconditionnel}.


\bigskip

Our estimates rely on a very simple observation. To state it and for further use in this paper, we give the following definition.

\begin{defi} A differential operator $ B $ is of 'div-grad' type if it is of the form
\begin{eqnarray}
 B = \sum_{j,k = 1}^d D_j \left( b_{jk}(x) D_k \right),  \label{operateurausensdesformes}
\end{eqnarray}  with coefficients such that $ b_{jk} \in S_{\rm dil} (\Ra^d) $. As usual, we have set $ D_j = \frac{1}{i} \frac{\partial}{\partial x_j} $.
\end{defi}

The first ingredient of the proof of Theorem \ref{theoremeprincipal} is the following trivial remark.
\begin{lemm} If $ B  $ is of div-grad type then $ [A,B] $ is of div-grad type. More precisely, if 
$$ B = \sum D_j (b_{jk}(x) D_k) , $$
then
\begin{eqnarray}
 i [B,A] = \sum_{jk} D_j \left( 2 b_{jk}(x) - (x \cdot \nabla b_{jk})(x) \right) D_k . \label{formeducommutateurtypique}
\end{eqnarray}
\end{lemm}
We omit the result which follows from an elementary computation (see also (\ref{formuletrivialeimportante}) below). Note that the formal computations are justified by the assumption that the coefficients $ b_{jk} $ are smooth. 

The second ingredient is the Mourre theory (see for instance \cite{Mourre}). 
Basically, the Mourre theory allows to derive a priori bounds on the solutions to
$$ (P - z) u = f , $$
(or more general Schr\"odinger operators), by exploiting a positive commutator estimate of the form
$$ \chi (P) i [P, A] \chi (P ) \geq c \chi^2 (P) , $$
with $ c > 0$ and $ \chi \in C_0^{\infty}(\Ra) $ real valued and equal to $1$ in a neighborhood of $ \mbox{Re}(z) $. For operators of div-grad type as in this paper, such estimates hold only if $ \chi $ is supported in $ \Ra^+ $, ie away from the $ 0 $ threshold. This is due to the fact that $ i [P,A] $ is close to $ 2 P $ (at least for globally small perturbations or near infinity), so one can essentially bound from below the (spectrally localized) commutator by
$ 2 \chi (P) P \chi (P) $. The latter is only positive definite (on the range of $ \chi (P) $) if $ \chi $ is supported in $ \Ra^+ $ and one then has $ || P^{1/2} \chi (P) v ||_{L^2} \approx || \chi (P) v ||_{L^2} $ by the spectral theorem. If $ 0 $ belongs to the support of $ \chi $, we loose this equivalence.  Rather than getting lower bounds by $ L^2 $ norms, we shall use the weaker observation that (in the simple case of small perturbations)
$$ (i[P , A ] v , v ) \geq || \nabla v ||_{L^2}^2 \gtrsim || v ||_{L^{2^*}}^2 $$
by the homogeneous Sobolev embedding
\begin{eqnarray}
 || v ||_{L^{2^*}} \leq C \big| \big| |D|v \big| \big|_{L^2} . \label{Sobolevembedding}
\end{eqnarray}
In other words, we keep the $ P^{1/2} $ factor to bound $ 2( \chi (P) P \chi (P) v,v)  $ from below by $ || P^{1/2} \chi (P) v ||^2 $.
By combining this remark with techniques due basically to Mourre, we shall derive (weighted) $ L^{2_*} \rightarrow L^{2^*} $
bounds for the resolvent of $ P $.

\section{Properties of the generator of dilations}
\setcounter{equation}{0}
In this section we collect some elementary formulas for the generator of dilations (\ref{generateur}) and its resolvent. For further purposes, it will be convenient to consider its semiclassical version, ie $ h A $ with $ 0 < h < 1 $.
All the properties will follow from the usual formula
\begin{eqnarray}
(hA-z)^{-1} = \frac{1}{i} \int_0^{\pm \infty} e^{-itz} e^{ithA} dt, \qquad \pm \mbox{Im}(z) < 0 , \label{semigrouperesolvente}
\end{eqnarray}
combined with the explicit form of the unitary group (\ref{groupe}).

Observe first that, since
\begin{eqnarray}
|| e^{i t h A} \varphi ||_{L^p} = e^{ht \left(\frac{d}{2} - \frac{d}{p}\right)}|| \varphi ||_{L^p}
\end{eqnarray}
for $ p \in [1,\infty] $ and, for instance, $ \varphi \in {\mathcal S}(\Ra^d) $, the formula (\ref{semigrouperesolvente}) implies that
\begin{eqnarray}
 || (h A - z)^{-1} \varphi ||_{L^p} \leq \frac{1}{|{\rm Im}(z)| - h |\frac{d}{2}-\frac{d}{p}|} || \varphi ||_{L^p} ,
\end{eqnarray}
provided that $ |\mbox{Im}(z) | > h |\frac{d}{2}-\frac{d}{p}| $. For the applications in this paper, this will be always the case since $ z $ will be close to $ \pm i$ and $h$ will be small.

Next, if $ \rho $ is a measurable function of polynomial growth, one readily checks that
\begin{eqnarray}
e^{ithA} \rho (D) e^{-ithA} & = & \rho (e^{-th}D) , \\
e^{ithA} \rho (x) e^{-ithA} & = & \rho (e^{ht}x) . \label{conjugaisonfonction}
\end{eqnarray}
Also, if $ \rho $ is $ C^1 $ with gradient of polynomial growth, we have
\begin{eqnarray}
i [\rho(D),A] & = & \left( \xi \cdot \nabla_{\xi} \rho \right) (D), \label{commutateurxi} \\
i [\rho(x),A] & = & - \left( x \cdot \nabla_x \rho \right) (x) . \label{commutateurx}
\end{eqnarray}
In the special case where $ \rho = \rho_s $ is homogeneous of real degree $ s \geq 0 $, we have
\begin{eqnarray}
e^{ithA} \rho_s (D) e^{-ithA} & = & e^{-sth} \rho_s (D) , \label{Fouriermultiplier} 
\end{eqnarray}
from which one easily deduce that  
\begin{eqnarray}
(hA - z)^{-1} \rho_s (D) = \rho_s (D) (hA - z + i h s)^{-1}, \qquad |\mbox{Im}(z)| > hs , \label{pourtraverserderivees}
\end{eqnarray}
using (\ref{semigrouperesolvente}). 

\bigskip

Finally, we consider the action on differential operators. If $ B = \sum_{jk} D_j \left( b_{jk}(x) D_k \right)  $ is of div-grad type, (\ref{conjugaisonfonction}) and (\ref{Fouriermultiplier})  readily imply that
\begin{eqnarray}
e^{ithA} B e^{-ithA} = e^{-2ht} \sum_{jk} D_j \left( b_{jk}(e^{ht} x) D_k \right) . \label{formuletrivialeimportante}
\end{eqnarray}
Operators of this form will be of great importance in this paper. Let us record the following simple property.
\begin{prop} \label{simplecalculation} Let $ b \in S_{\rm dil} (\Ra^d) $ and set $$ b_{(\tau)}(x) = b (e^{\tau}x) , $$
ie $ b_{(\tau)} = e^{i\tau A} b e^{-i\tau A} $ as multiplication operators. Then, for all $ k , n \in \Na $,
$$ \partial_{\tau}^k (x \cdot \nabla)^n \left( b_{(\tau)} \right) = ((x \cdot \nabla)^{k+n} b)_{(\tau)}  . $$
In particular, for all $ N $,
\begin{eqnarray}
 || b ||_{N,{\rm dil}} = || b_{(\tau)} ||_{N,{\rm dil}} . \label{stabiliteparconjugaison}
\end{eqnarray}
\end{prop}

\noindent {\it Proof.} A straightforward calculation which we omit. \finpreuve

\bigskip

For further purposes, it will be convenient to use the following definition.
\begin{defi}[Admissible operators] \label{definitionadmissible} Let $ m \in \Na $. We say that a family $ ( b_{\tau} )_{\tau \in \Ra} $ is $ m $-admissible in $   S_{\rm dil} (\Ra^d) $ if, for all integers $ k, n $
$$  || \partial_{\tau}^k (x \cdot \nabla)^n b_{\tau}  ||_{L^{\infty}} \leq C_{kn}e^{m|\tau|} . $$
A family of differential operators $ ( B_{\tau} )_{\tau \in \Ra} $ is $ m $-admissible if 
\begin{eqnarray}
 B_{\tau} =  \sum_{j,k=1}^d D_j(b_{jk,\tau}(x)  D_k) , \nonumber \label{formeadmissible}
\end{eqnarray}
with $ (b_{jk,\tau})_{\tau \in \Ra} $ $ m $-admissible families in $ S_{\rm dil}(\Ra^d) $.
\end{defi}

\bigskip

\noindent {\bf Example.} With the notation of Proposition \ref{simplecalculation}, $ b_{\tau}^{\pm} := e^{ \pm 2 \tau} b_{(\tau)} $ are two $ 2$-admissible families in $ S_{\rm dil} (\Ra^d) $. 


\bigskip

\begin{prop} \label{pourlesintegrationsparparties} Let $ (B_{\tau})_{\tau \in \Ra} $ be a $m$-admissible family of differential operators. Then, if $ w : [0,1] \rightarrow \Ca $ is continuous, the operators
\begin{eqnarray}
 \frac{d}{d\tau} B_{\tau} , \qquad e^{i \tau A} B_{\tau} e^{-i\tau A} \qquad \mbox{and} \qquad \int_0^1 w(s) B_{s \tau} ds , \nonumber
 \end{eqnarray}
 are respectively $ m $, $ m +2 $ and $ m $-admissible.
\end{prop}

In this proposition, the derivative $ \frac{d}{d \tau} $ (resp. integration) mean that one considers the operator with coefficients differentiated (resp. integrated) with respect to $ \tau $.

\bigskip

\noindent {\it Proof.} The case of $ (d/d\tau)B_{\tau} $ is obvious. For the second operator, the result follows from (\ref{formuletrivialeimportante}) (with $ th = \tau $) and the fact that $ m $-admissible coefficients are stable by conjugation by $ e^{i \tau A} $ which is due to Proposition \ref{simplecalculation}.  The last case is simply a consequence of the fact that $ \int_0^1 |w(s)|s^k e^{m|s\tau|} ds \lesssim e^{m|\tau|} $, for all non negative integer $k$. \finpreuve

\section{A representation formula for the commutator}
\setcounter{equation}{0}
As indicated in the introduction, we shall use the commutator techniques of Mourre to get lower bounds. It will be convenient to use the recent energy estimates approach proposed by G\'erard \cite{Gerard}. The purpose of the present section is to compute relatively explicitly the relevant commutator.

In the sequel we denote by $ F $ the bounded function
$$ F (\lambda) = \arctan (\lambda), \qquad \lambda \in \Ra , $$
whose final interest will be that it is positive (or negative) up to an additive constant and has a positive derivative.

We also introduce
\begin{eqnarray}
 P_{\tau} = e^{- i \tau A} i [ P, A ] e^{i \tau A} , \label{pourlecalculdesderivees}
\end{eqnarray}
and standardly denote
\begin{eqnarray}
 \left( i[P,F(hA)]u_1,u_2 \right) = (i F (h A) u_1 , P u_2) - (i P u_1 , F (h A) u_2) . \label{formsense}
\end{eqnarray}
The purpose of this section is to prove a representation formula for this commutator. Rather than using the Helffer-Sj\"ostrand formula
as in \cite{GoleniaJecko}, we use here a functional calculus based on Fourier transform which is more convenient since we have an explicit formula for the unitary group $ e^{itA} $.

\begin{prop} \label{formulederepresentation} For all $ u_1 , u_2 \in {\mathcal S}(\Ra^d) $ and all $ 0 < h < 1 $, we have
\begin{eqnarray}
 \left( i[P,F(hA)]u_1,u_2 \right) = \frac{h}{2} \int_{\Ra} e^{-|t|}  \left( \frac{1}{t} \int_0^t \left( e^{ithA} P_{sh} u_1 , u_2 \right) ds \right) dt . 
 \label{formuleexacte}
\end{eqnarray}
\end{prop}

In the spirit of \cite{Gerard}, we use a semiclassical parameter $h$  thanks to which  the derivation of a positive estimate will be fairly transparent.

\bigskip

The rest of the section is devoted to the proof of this proposition. Recall first that
$$ \arctan (\lambda) = \int_0^{+\infty} \frac{\sin(t\lambda)}{t} e^{-t} dt , $$
which we are going to approximate by 
$$ F_{\nu} (\lambda) = \int_0^{+\infty} \sin(t\lambda) \frac{t}{t^2 + \nu^2} e^{-t} dt = \frac{1}{2i} \int_{\Ra} e^{i t \lambda} \frac{t}{t^2 + \nu^2} e^{-|t|} dt , $$
with $ \nu > 0 $. For future reference, we record here the following lemma.
\begin{lemm} \label{convergenceforte} There exists $ C > 0 $ such that
\begin{eqnarray}
|F_{\nu}(\lambda)| \leq C |\lambda|, \qquad \nu > 0, \ \lambda \in \Ra .
\end{eqnarray}
Furthermore, for all $ \lambda \in \Ra $,
\begin{eqnarray}
 F_{\nu}(\lambda) \rightarrow F (\lambda), \qquad \nu \rightarrow 0 . 
 \end{eqnarray}
\end{lemm}

We omit the very simple proof.

\bigskip

\begin{lemm} For all $ v , w \in L^2 (\Ra^d) $, all $ \nu > 0 $ and all $ h > 0 $, we have
\begin{eqnarray}
 (F_{\nu}(hA)v,w) = \frac{i}{2} \int_{\Ra} \frac{t e^{-|t|}}{t^2 + \nu^2}  \left( e^{ithA}v,w \right) dt . \label{avecParseval}
\end{eqnarray}
\end{lemm}

\noindent {\it Proof.} If $ ( E^{hA}_{\lambda} )_{\lambda \in \Ra} $ denotes the spectral resolution of $hA$, we have by definition
$$ (F_{\nu}(hA)v,w) = \int_{\Ra} F_{\nu} (\lambda) d \left( E_{\lambda}^{hA} v , w \right) , $$
and then by Parseval's identity
\begin{eqnarray*}
 (F_{\nu}(hA)v,w)  & = & \frac{1}{2 \pi} \int_{\Ra} \widehat{F}_{\nu}(t) (e^{-ithA}v,w) dt , \\
 & = & \frac{i}{2} \int_{\Ra} \frac{t e^{-|t|}}{t^2 + \nu^2}  \left( e^{ithA}v,w \right) dt .
\end{eqnarray*}
This identity can be justified by a standard density argument, assuming first that $ v $ and $w$ are spectrally localized (ie of the form $ \chi(A)v $, $ \chi (A) w$ with $ \chi \in C_0^{\infty} $) and approximating (for fixed $ \nu $) $ F_{\nu} $ by Schwartz functions by adding a cutoff vanishing close to $ t=  0 $ in the definition of $ F_{\nu} $. These Schwartz functions converge pointwise to $ F_{\nu} $ with uniform bound of order $ C|\lambda| $ which is harmless if we consider spectrally localized $v$ and $w$. Their fourier transform  converge $ dt$ almost everywhere (pointwise on $ \Ra_t \setminus 0$) to $ \widehat{F}_{\nu} $ with uniform bound by $ C |t|e^{-|t|} $ and the result follows then easily.
 \finpreuve

\bigskip

Since $ F_{\nu} $ is real valued, we have $ (F_{\nu}(hA)v,w)= (v,F_{\nu}(hA)w) $ and thus
\begin{eqnarray}
 (v,F_{\nu}(hA)w) = \frac{i}{2} \int_{\Ra} \frac{t e^{-|t|}}{t^2 + \nu^2} \left( v , e^{-ithA} w \right) dt . \nonumber
\end{eqnarray}
From the latter identity and (\ref{avecParseval}), we deduce that
\begin{eqnarray}
\left( i [P , F_{\nu}(hA)] u_1 , u_2 \right) = \frac{1}{2} \int \frac{te^{-|t|}}{t^2 + \nu^2} \left( (e^{ithA}u_1,P u_2) - (Pu_1,e^{-ithA}u_2) \right) dt ,
\label{pourpasseralalimite}
\end{eqnarray}
where the commutator in the left hand side is understood in the  sense of (\ref{formsense}) (ie the form sense).

\begin{lemm} For all $ t \in \Ra $, $ h > 0 $ and $ u_1 , u_2 \in {\mathcal S}(\Ra^d) $,
\begin{eqnarray}
(e^{ithA}u_1,P u_2) - (Pu_1,e^{-ithA}u_2) = h \int_0^t \left( e^{ithA} P_{sh} u_1 , u_2 \right) ds . \label{pouralgebre}
\end{eqnarray}
In addition, for each pair $ u_1, u_2 $, there is a constant $ C $ such that
\begin{eqnarray}
\left| (e^{ithA}u_1,P u_2) - (Pu_1,e^{-ithA}u_2) \right| \leq C|t| e^{ h|t|}, \qquad t \in \Ra . \label{pourconvergencedominee}
\end{eqnarray}
\end{lemm}

\noindent {\it Proof.} The formula (\ref{pouralgebre}) is equivalent to the same one with $ u_1 $ replaced by $ e^{-ithA} u_1 $ and the corresponding identity is then a consequence of Duhamel's formula, ie is obtained by checking that the derivatives of both sides coincide, using (\ref{formuletrivialeimportante}). To get (\ref{pourconvergencedominee}), we use (\ref{pouralgebre}) and observe that, since the coefficients of $ P_{sh} $ are of order $ e^{2sh} $ (see (\ref{formuletrivialeimportante}) and (\ref{pourlecalculdesderivees})), we have
\begin{eqnarray*}
| \left( e^{ithA} P_{sh} u_1 , u_2 \right) | & \leq & C e^{2sh} || \nabla  u_1 ||_{L^2} || \nabla e^{-ithA }  u_2 ||_{L^2} \\
 & \leq & e^{(2s-t)h} || \nabla u_1 ||_{L^2} || \nabla u_2 ||_{L^2}
\end{eqnarray*}
where $ |2s-t| \leq |t| $ since $s$ is between $0$ and $t$. The conclusion follows easily.
\finpreuve

\bigskip

\noindent {\it Proof of Proposition \ref{formulederepresentation}.} By Lemma \ref{convergenceforte} and the Spectral Theorem, we have 
$$ F_{\nu}(hA) u_j \rightharpoonup F (h A) u_j, \qquad \nu \rightarrow 0 , \ j = 1,2. $$
 Thus the left hand side of (\ref{formuleexacte}) is the limit as $ \nu \rightarrow 0 $ of the left hand side of (\ref{pourpasseralalimite}).
To compute the limit  of the right hand side of (\ref{pourpasseralalimite}), we simply insert (\ref{pouralgebre}) therein and then let $ \nu \rightarrow 0 $ by dominated convergence using (\ref{pourconvergencedominee}) and the fact that $ h < 1 $. The limit is clearly the right hand side of (\ref{formuleexacte}) and this completes the proof. \finpreuve

\section{Semiclassical expansion of the commutator}
\setcounter{equation}{0}
In this section, we establish the first order asymptotic expansion in $h$ of (\ref{formuleexacte}). To state this result, we introduce the following notation.
Write first 
\begin{eqnarray}
 P_{sh} = P_0 + sh Q_{sh} , \label{decomposition}
\end{eqnarray}
with
\begin{eqnarray}
 Q_s = \int_0^1 \frac{d}{d \tau} P_{\tau|_{\tau = \sigma s } } d \sigma . \label{pourlecalculdesderivees2}
\end{eqnarray}
Write next
$$ \frac{1}{t} \int_0^t h s Q_{sh} ds = t h \int_0^1  s  Q_{t s h} d s ,  $$
and set 
\begin{eqnarray}
 B_{\tau} : = \tau \int_0^1 s Q_{s \tau} ds . \label{formeRtau}
\end{eqnarray}
Notice that $ (P_{\tau})_{\tau \in \Ra} $ given by (\ref{pourlecalculdesderivees}) is a $ 2 $-admissible family of differential operators (see Definition \ref{definitionadmissible})  hence so are $ (Q_{\tau})_{\tau \in \Ra} $ and $ (B_{\tau})_{\tau \in \Ra} $ by Proposition \ref{pourlesintegrationsparparties}. 

Observe that
\begin{eqnarray}
 \frac{h}{2} \int e^{-|t|} ( e^{ithA} P_0 u_1 , u _2 ) dt = h \left( P_0 u_1 , (h^2 A^2 + 1)^{-1} u_2 \right) , \label{exactformula}
\end{eqnarray}
as follows easily from the spectral theorem and the Fourier transform 
$$ \frac{1}{1 + \lambda^2} = \frac{1}{2} \int_{\Ra} e^{-it \lambda} e^{-|t|} dt . $$
It can also be seen as a consequence of (\ref{semigrouperesolvente}). Define
$$ {\mathcal A}_{h,H} (u_1,u_2) := \left( P_0 (hA+i)^{-1}u_1 , (h A + i)^{-1} u_2 \right) , $$
and
$$ {\mathcal B}_{H,h} (u_1,u_2) = \frac{1}{h} \left\{ ( P_0 u_1 , (h^2 A^2 + 1)^{-1} u_2 ) - (P_0 (hA+i)^{-1}u_1,(hA+i)^{-1} u_2) \right\} , $$
so that
$$ h \left( P_0 u_1 , (h^2 A^2 + 1)^{-1} u_2 \right) = h {\mathcal A}_{h,H}(u_1,u_2) + h^2 {\mathcal B}_{h,H} (u_1,u_2) . $$
If we finally set
$$ ( {\mathcal C}_h u_1 , u_2 ) :=  \frac{1}{2} \int_{\Ra} e^{-|t|} t (e^{i t h A} B_{th} u_1 , u_2 ) dt , $$
we have
\begin{eqnarray}
\left( i[P,F(hA)]u_1,u_2 \right)  =  h {\mathcal A}_{h,H}(u_1,u_2) + h^2 {\mathcal B}_{h,H} ( u_1 , u_2) + h^2 {\mathcal C}_{h,H} (u_1,u_2) .  \label{reformulation}
\end{eqnarray}
The purpose of this section is thus to estimate $ {\mathcal B}_{h,H} $ and $ {\mathcal C}_{h,H} $.

\begin{prop} \label{fonctionnel1} There exists $C $ such that for all $ 0 < h < 1 $ and all $ H $, 
\begin{eqnarray}
\left| {\mathcal B}_{H,h} (u_1,u_2) \right| \leq C  (1+||H||_{2,{\rm dil}}) \big| \big| |D| (hA+i)^{-1} u_1 \big| \big|_{L^2} \big| \big| |D| (hA+i)^{-1} u_2 \big| \big|_{L^2} .
\end{eqnarray}
\end{prop}

\noindent {\it Proof.} By the resolvent identity
\begin{eqnarray}
(hA+i+ih)^{-1} = (hA + i )^{-1} - ih (hA+i)^{-1}(hA+i+ih)^{-1}
\end{eqnarray}
and (\ref{pourtraverserderivees}), we have
\begin{eqnarray}
 (hA+i)^{-1} D_j = D_j \left( 1 - i h (hA+i+ih)^{-1} \right) (hA+i)^{-1}  . \label{commutateur0}
\end{eqnarray}
Next, we observe that
\begin{eqnarray}
 [(hA+i)^{-1}, G_{jk}] = - \frac{h}{i} (hA+i)^{-1} (x \cdot \nabla H_{jk}) (hA+i)^{-1} , \label{commutateur1}
\end{eqnarray}
and finally that
we also have
\begin{eqnarray}
(hA+i)^{-1}D_k = \left(1 - ih (hA+i)^{-1} \right) D_k (hA+i)^{-1} , \label{commutateur2}
\end{eqnarray}
 since
\begin{eqnarray}
[(hA+i)^{-1} , D_k  ] = - h (hA+i)^{-1}[A,D_k] (hA+i)^{-1} =  i h (hA+i)^{-1} D_k (hA+i)^{-1} .  \nonumber
\end{eqnarray}
From (\ref{commutateur0}), (\ref{commutateur1}) and (\ref{commutateur2}), we see that
$$ [ (hA+i)^{-1} , P_0 ] = \sum_{jk} D_j B_{jk}(h) D_k (hA+i)^{-1} , $$
with
$$ || B_{jk}(h) ||_{L^2 \rightarrow L^2} \lesssim  h \left( 1 + ||H||_{2,{\rm dil}} \right) . $$
The result follows. \finpreuve

\bigskip







\begin{prop} \label{integrationresolvente} For all $ 0 < h_0 < 1 / 4 $, there exists $ C > 0 $ such that
$$ | {\mathcal C}_{h,H} (u_1,u_2) |\leq C (1+||H||_{4,{\rm dil}}) \big|\big| |D|(hA+i)^{-1} u_1||_{L^2} \big|\big| |D|(hA+i)^{-1} u_2 \big|\big|_{L^2} $$
for all $ u_1, u_2 \in {\mathcal S}(\Ra^d) $, all $ 0 < h < h_0 $ and all $ H $.
\end{prop}



\noindent {\it Proof.} It simply relies on integrations by parts. Indeed, since
\begin{eqnarray}
e^{-ithA} u_2 = i e^{-ithA} (hA+i)^{-1} u_2 + i \frac{d}{dt} e^{-ithA} (hA+i)^{-1} u_2 \label{areutiliserpouru1}
\end{eqnarray}
we can write
$$ {\mathcal C}_{h,H} (u_1,u_2) = i {\mathcal C}_{h,H}(u_1, (hA+i)^{-1}u_2) + \frac{i}{2} \int t e^{-|t|} \left( B_{ht} u_1 , \frac{d}{dt} e^{-ithA} (hA+i)^{-1} u_2 \right) dt , $$
where the second term in the right hand side reads
$$ - \frac{i}{2} \int  e^{-|t|} \left( \{ t  h B_{ht}^{\prime} + (1-|t|)B_{ht} \}  u_1 ,  e^{-ithA} (hA+i)^{-1} u_2 \right) dt , $$
if $ B^{\prime}_{\tau} = (d/d\tau) B_{\tau} $. Recall that $ (B^{\prime}_{\tau})_{\tau} $ is still a $2$-admissible family of operators so that
$$ \tilde{B}_{\tau}:=e^{i \tau A} B_{\tau} e^{-i \tau A} \qquad \hat{B}_{\tau}:= e^{i \tau A} B_{\tau}^{\prime} e^{-i \tau A} , $$
define $4 $-admissible families of operators by Proposition \ref{pourlesintegrationsparparties}. Then, 
using again (\ref{areutiliserpouru1}) with $ u_1 $ instead of $ u_2 $ and integrating by parts (remark that
the functions $ e^{-|t|} $ and $ (1-|t|) $ are not $ C^1 $ at $ t = 0 $ but  are continuous and therefore there are no boundary terms), we obtain a sum of integrals of the form
$$  \int_0^{\pm \infty} w_{\pm}(t)  e^{-|t|} \left( e^{ithA} C_{ht}^{\pm} (hA+i)^{-1}  u_1 ,  e^{-ithA} (hA+i)^{-1} u_2 \right) dt $$
with $ w_{\pm} $ polynomial and $ ( C_{\tau}^{\pm} )_{\tau \in \Ra} $ $ 4 $-admissible families of operators whose coefficients are bounded in $ L^{\infty}(\Ra^d) $ by $ e^{4 |\tau|} || H ||_{4,{\rm dil}} $. The result follows.  \finpreuve

\section{Proofs of the results}
\setcounter{equation}{0}
\subsection{Proof of Theorem  \ref{theoremeprincipal} }
Assume  that $ \mbox{Im}(z) > 0 $. The estimates for $ \mbox{Im}(z) < 0 $ are obtained by taking the adjoint. We recall that $ F (\lambda) = \arctan (\lambda) $. As in \cite{Gerard}, we observe that
\begin{eqnarray}
2 \mbox{Im} \left( \left( F (hA)-\frac{\pi}{2} \right) u , (P-z) u \right) & = & 2 \mbox{Im}(  F (hA) u , P u) - 2 \left( \mbox{Im}(z) \left( F(hA) - \frac{\pi}{2} \right) u, u \right)  \nonumber \\
 & = & (i[P,F(hA)]u,u) - 2 \left( \mbox{Im}(z) \left( F(hA) - \frac{\pi}{2} \right) u, u \right) \nonumber \\
 & \geq & (i [P,F(hA)]u,u) .
\end{eqnarray}
By Propositions \ref{formulederepresentation}, \ref{fonctionnel1} and \ref{integrationresolvente}, we have
\begin{eqnarray}
(i [P,F(hA)]u,u) & \geq & h (P_0 (h A + i)^{-1}u,(hA+i)^{-1}u) - C h^2 || |D| (h A + i)^{-1} u ||_{L^2}^2 \label{premiereligne} \\
& \geq & \frac{\varepsilon}{2} h \big| \big| |D| (hA+i)^{-1} u \big| \big|_{L^2}^2 , \nonumber
\end{eqnarray}
by taking $ h $ small enough so that $ C h \leq \varepsilon / 2 $. Notice that the constant $ C $ in (\ref{premiereligne}) is of order $ 1 + ||H||_{4,{\rm dil}} $ so that we may choose $h^{-1}$ of order $ (1 + ||H||_{4,{\rm dil}}) $.

On the other hand, we may write 
$$  \left( \left( F (hA)-\frac{\pi}{2} \right) u , (P-z) u \right) = \left(|D| \left( F (hA)-\frac{\pi}{2} \right) (hA+i)^{-1} u ,|D|^{-1}(hA-i) (P-z) u \right)  . $$
Thus, once we have proved Proposition \ref{propositionFhA} below,  we shall get the estimate
$$ \big| \big| |D|(hA+i)^{-1} u \big| \big|_{L^2} \leq \frac{C}{h} \big| \big| |D|^{-1}(hA-i) (P-z) u \big| \big|_{L^2} $$
which gives (\ref{inegaliteprincipale}).  
\begin{prop} \label{propositionFhA} For all $ 0 < h_0 < 1 $, there exists $ C > 0 $ such that
$$ \big| \big| |D| F(hA) (hA+i)^{-1} u \big| \big|_{L^2} \leq C \big| \big| |D| (hA+i)^{-1} u  \big| \big|_{L^2} , $$
for all $ u \in {\mathcal S}(\Ra^d) $ and $ 0 < h \leq h_0 $.
\end{prop}

\noindent {\it Proof.} Since we have
$$ \big| \big|  F(hA) |D| (hA+i)^{-1} u \big| \big|_{L^2} \leq ||F||_{\infty} \big| \big| |D| (hA+i)^{-1} u  \big| \big|_{L^2} , $$
the result is clearly equivalent to an estimate on the commutator $ [|D|,F(hA)] $. The latter can be computed explicitly using the same argument as for
  Proposition \ref{formulederepresentation}. We obtain
$$ \left( i [|D|,F(hA)]u_1,u_2 \right) = \frac{h}{2} \int_{\Ra} e^{-|t|} \left( \frac{1}{t} \int_0^t e^{sh}( e^{ithA} |D| u_1 , u_2 ) ds \right) dt, \qquad
u_1,u_2 \in {\mathcal S}(\Ra^d) , $$
since, $ e^{-ishA} i [|D|,A] e^{ishA} = e^{sh}|D| $. This implies that
$$ \big| ([|D|,F(hA)]u_1,u_2 ) \big| \leq \frac{h}{2} \int e^{-(1-h)|t|} dt \big| \big| |D| u_1 \big| \big|_{L^2} ||u_2||_{L^2} , $$
ie that $ || [|D|,F(hA)]u_1 ||_{L^2} \lesssim (1-h)^{-1} \big| \big| |D| u_1 \big| \big|_{L^2} $. The result then follows clearly.
 \finpreuve
 
\subsection{Proof of Corollary \ref{corollaireTao}}
Using the homogeneous Sobolev imbedding (\ref{Sobolevembedding}), we have, for any $ f \in L^2$ 
\begin{eqnarray}
 \big| \big| (hA+i)^{-1} (P-z)^{-1}  f \big| \big|_{L^{2^*}} \leq C \big| \big| |D|
 (hA+i)^{-1} (P-z)^{-1}  f \big| \big|_{L^{2}} . \label{acombiner}
\end{eqnarray}
 Then, by  choosing $f  = (hA-i)^{-1} g $ with $ g \in L^2 \cap L^{2_*} $, we have 
 \begin{eqnarray*}
  \big| \big| |D|
 (hA+i)^{-1} (P-z)^{-1}  f \big| \big|_{L^{2}} & = & \sup_{||u||_{L^2} = 1} \big| \left( |D|
 (hA+i)^{-1} (P-z)^{-1}  f ,u \right) \big| \\
 & = & \sup_{||u||_{L^2} = 1} \big| \left( g,
 (hA+i)^{-1} (P-\bar{z})^{-1} (hA-i)^{-1} |D|  u \right) \big| \\
 & \leq & \sup_{ ||u||_{L^2} = 1 } ||g||_{L^{2_*}} \big| \big| (hA+i)^{-1}(P-\bar{z})^{-1} (hA-i)^{-1} |D|  u \big| \big|_{L^{2^*}} \\
 & \leq &  \big| \big| |D|(hA+i)^{-1}(P-\bar{z})^{-1} (hA-i)^{-1} |D|  \big| \big|_{L^2 \rightarrow L^2} ||g||_{L^{2_*}}
\end{eqnarray*}
which combined with (\ref{acombiner}) completes the proof. \finpreuve
\subsection{Proof of Corollary \ref{longueportee}}
By H\"older's inequality, 
$$   ||\scal{x}^{-1 - \epsilon} u ||_{L^2} \lesssim || u ||_{L^{2^*}}, \qquad || \scal{x}^{-1- \epsilon} v ||_{L^{2_*}} \lesssim ||v||_{L^2} .  $$
Choose next $ \chi \in C_0^{\infty}(\Ra) $ which is equal to $ 1 $ near $ [0,1] $. It is classical that
\begin{eqnarray}
(1 - \chi^2)(P) (P-z)^{-1} : L^{2_*} \rightarrow L^{2^*}
\end{eqnarray}
by Sobolev embeddings, with norm uniformly bounded for $ |z| \leq 1 $. This follows for instance from the fact that the $L^2$ bounded operator $ (1 - \chi^2)(P) (P-z)^{-1} $ is a pseudo-differential operator of order $ -2 $. It is therefore sufficent to show that
$$ \scal{x}^{-1} \chi (P) (P-z)^{-1} \chi (P) \scal{x}^{-1} : L^{2_*} \rightarrow L^{2^*} $$
is bounded uniformly with respect to $ |z|<1 $, $  z \notin \Ra $. To get the latter, we simply write
$$ \chi (P) \scal{x}^{-1} = (hA-i)^{-1} (hA-i) \chi (P) \scal{x}^{-1} $$
and use the fact that
$ \chi (P) \scal{x}^{-1} $  and $ A \chi (P) \scal{x}^{-1}  $ are bounded on $ L^{p} $ for all $p$, which follows from the fact that these operators are pseudo-differential operators of order $ - \infty $ (see for instance \cite{stBoTz2} for more details on such properties).\finpreuve

\subsection{Local compactness of the Riesz transform} \label{soussectionRiesz}
In this subsection we prove of  property of the Riesz transform which we shall use in the proof of Theorem \ref{theoremeconditionnel}.
 We first recall the definition of the Riesz transform. Since $ P \geq 0 $ is self-adjoint, the spectral theorem and (\ref{bornesellipticite}) give
\begin{eqnarray}
  (Pu,u) = || P^{1/2} u ||_{L^2}^2 \approx || \nabla u ||_{L^2}^2 \approx \big| \big| |D| u \big| \big|_{L^2}^2 , \label{controlederivationbassefrequence}
\end{eqnarray}
where $ \approx $ stands for the equivalence of norms. This implies that $ P^{1/2} $ is an isomorphism from the homogeneous Sobolev space $ \dot{H}^1(\Ra^d) $ onto $ L^2 (\Ra^d) $. Denoting the inverse by $ P^{-1/2} $, the operators 
\begin{eqnarray}
 R (j) = \partial_{x_j} P^{-1/2} , \label{transformeedeRiesz}
\end{eqnarray} 
  are then well defined on $ L^2 (\Ra^d) $ for all $j$. They are the components of the well known  Riesz transform $ \nabla P^{-1/2} $. To define more explicitely $ R(j) $, we can use the following integral representation (see for instance \cite{Auscher}). For each $n \geq 1$, we consider
$$ R_n (j) := \pi^{-1/2} \partial_{x_j} \int_{1/n}^n e^{-t P} \frac{dt}{\sqrt{t}} , $$
where the integral converges in the strong sense. It is not hard to check that $ R_n (j) $ is bounded using that $ e^{-tP} $ maps $ L^2 $ in $ \cap_s H^s $ for all $t>0$. Let us briefly recall why $ R_n (j) $ converges strongly as $n \rightarrow \infty$ (for this purpose we could actually consider  lower and upper bounds in the integral defining   $ R_n (j) $ going independently to $ 0 $ and $ \infty $ respectively,  but this is irrelevant for our purpose). Using (\ref{controlederivationbassefrequence}), we see that
\begin{eqnarray}
 || R_n (j) u ||_{L^2} \leq C \left| \left| \sqrt{P} \int_{1/n}^n e^{-tP} \frac{dt}{\sqrt{t}} u \right| \right|_{L^2} =  C|| f_n (P) u ||_{L^2} , \label{bornespectrale}
\end{eqnarray}
with
$$ f_n (\lambda) = \int_{1/n}^n \lambda^{1/2} e^{-t \lambda} \frac{dt}{\sqrt{t}} = \int_{\lambda/n}^{\lambda n} e^{- \tau} \frac{d \tau}{\sqrt{\tau}} . $$
Since $f_n$ is uniformly bounded with respect to $ n \geq 1 $ and $ \lambda \geq 0 $, (\ref{bornespectrale}) and the spectral theorem show that $ || R_n (j) ||_{L^2 \rightarrow L^2} \leq C $ for all $n$. Therefore, it sufficient to prove the strong convergence of $ R_n (j) $ on a dense subset.  For the latter, we observe that, since $ 0 $ is not an eigenvalue of $ P $, the spectral theorem shows that for all $ u \in L^2 $,
\begin{eqnarray}
 \chi_{[\epsilon , \epsilon^{-1}]}(P) u \rightarrow u, \qquad \epsilon \rightarrow 0 , \label{densite}
\end{eqnarray}
$ \chi_{[\epsilon , \epsilon^{-1}]} $ denoting the characteristic function of $ [\epsilon, \epsilon^{-1}] $. It is then easy to check that $ R_n (j) \chi_{[\epsilon , \epsilon^{-1}]}(P) $ converges in the strong sense a $n \rightarrow \infty $ for each $ \epsilon > 0 $ since the spectral projection on $ [\epsilon , \epsilon^{-1}] $ guarantees the exponential decay of $ e^{-t P} $ as well as the boundedness of $ \partial_{x_j} \chi_{[\epsilon , \epsilon^{-1}]}(P) $. By (\ref{densite}), functions of the form $ \chi_{[\epsilon , \epsilon^{-1}]}(P) u $ are dense $ L^2 $ so this completes the proof of the strong convergence of $ R_n (j) $. We may thus define
$$ R (j) = \pi^{-1/2} \partial_{x_j} \int_0^{\infty} e^{-tP} \frac{dt}{\sqrt{t}} := {\rm s-} \! \lim_{n \rightarrow \infty} R_n (j) , $$
which is a reasonable definition for $ \partial_{x_j} P^{-1/2}$ since one checks that
\begin{eqnarray}
 R (j) P^{1/2} u = \partial_{x_j} u,   \label{proprietealgebriqueRiesz}
\end{eqnarray}
for all $ u \in D (P) $. This is an elementary consequence of the spectral theorem and the Lebegue theorem since, for all $ \lambda > 0 $
$$ \pi^{-1/2} f_n (\lambda) \rightarrow 1, \qquad n \rightarrow \infty , $$
and since $ \{\lambda = 0 \} $ is negligible with respect to the spectral measure for $ 0 $ is not an eigenvalue of $P$. 
  This completes our definition of $ R (j) $. 
  
  The main purpose of the present subsection is to prove the following result.
\begin{prop} \label{conditionabstraite} Assume that $ d \geq 3 $. Then, for all $ \chi \in C_0^{\infty}(\Ra^d) $ and all $ \varphi \in C_0^{\infty}(\Ra) $,
\begin{eqnarray}
\chi (x) R (j) \varphi (P)  \ \ \mbox{is a compact operator on} \ L^2 (\Ra^d) , \nonumber
\end{eqnarray} 
for all $ j = 1, \ldots, d $.
\end{prop}
 
\noindent {\it Proof.} We split $ \pi^{1/2}R (j) $ into $ \partial_{x_j}\int_0^2 e^{-tP} dt/t^{1/2} + \partial_{x_j}\int_2^{\infty}  e^{-tP} dt/t^{1/2}$. It is clear that
$$ \chi (x) \partial_{x_j}\int_0^2 e^{-tP} \frac{dt}{\sqrt{t}} \chi (P) = \left( \chi (x) \partial_{x_j} \chi (P) \right) \int_0^2 e^{-tP} \frac{dt}{\sqrt{t}} $$
is compact since the bracket is compact and the integral defines a bounded operator on $ L^2 $. We then write contribution of the second term as
$$ \left( \chi (x) \partial_{x_j} e^{-P} \scal{x}^N \right) \int_2^{\infty} \scal{x}^{-N} e^{-(t-1)P} \frac{dt}{\sqrt{t}} , $$ 
with  $ N > 0 $ to be chosen below. Again the bracket is a compact operator. To see that the integral is bounded on $ L^2 $, we use the classical gaussian upper bounds for the kernel $ K (t,x,y) $ of $ e^{-tP} $ (see for instance \cite{Aronson,Davies}):  for some $ C,c > 0 $ we have,
\begin{eqnarray}
|K(t,x,y)| \leq \frac{C}{t^{d/2}} \exp \left(  \frac{c |x-y|^2}{t} \right) , \qquad x,y \in \Ra^d, \ t > 0, \nonumber
\end{eqnarray}
and thus
\begin{eqnarray}
||e^{-tP}||_{L^2 \rightarrow L^{\infty}} \lesssim t^{-d/4} .
\end{eqnarray}
Therefore, if $ N > d/2 $, 
$$ || t^{-1/2} \scal{x}^{-N} e^{-(t-1) P} u ||_{L^2} \lesssim t^{-\frac{1}{2}-\frac{d}{4}} ||u||_{L^2} , $$
which is integrable on $ [2,\infty) $ since $ \frac{1}{2} + \frac{d}{4} > 1 $. This completes the proof. \finpreuve

\subsection{Proof of Theorem \ref{theoremeconditionnel}}
We start by remarking that it is sufficient to show that, for some $ \lambda > 0 $ and $ h > 0 $ small enough, we have the bound
\begin{eqnarray}
 \big| \big| |D| (hA+i)^{-1} ( P - z )^{-1} (hA-i)^{-1} |D| \big| \big|_{L^2 \rightarrow L^2} \leq C, \qquad | \mbox{Re}(z) | < \lambda .
 \label{estimeereduite}
\end{eqnarray}
We will then obtain (\ref{bassefrequencelongueportee}) exactly as in Corollary \ref{longueportee}. We may even replace $ ( P - z )^{-1} $ in this estimate by $ (P-z)^{-1} \varphi_0 (P/\lambda) $, with $ \varphi_0 \in C_0^{\infty} (\Ra) $ such that $ \varphi_0 \equiv 1 $ near $ [-1,1] $, since the operator
$$ |D| (hA+i)^{-1} ( P - z )^{-1} (1-\varphi_0)(P/\lambda) (hA-i)^{-1} |D| $$
is easily seen to be bounded on $ L^2 $, uniformly with respect to $ z $ such that $ | \mbox{Re}(z) | < \lambda $.  It is therefore enough to consider $u$ of the form $ u =(P-z)^{-1} \varphi_0 (P /\lambda) f $ with $ f \in {\mathcal S}(\Ra^d) $ so that
\begin{eqnarray}
 u = \varphi (P / \lambda) u , \label{localisationspectrale}
\end{eqnarray}
for some $ \varphi \equiv 1 $ near $ \mbox{supp}(\varphi_0) $.

Independently, observe that, as in the proof of Theorem \ref{theoremeprincipal}, we have, for all $ u \in {\mathcal S}(\Ra^d) $,
\begin{eqnarray}
 (i [P,F(hA)]u,u) \geq \frac{h}{2} \left( P_0 (hA+i)^{-1} u , (hA+i)^{-1} u \right) - C h^2 \big| \big| |D| (hA+i)^{-1} u  \big| \big|_{L^2}^2 , 
 \label{commutateurnonelliptique}
\end{eqnarray}
but the  difference is now that $ P_0 $ is not necessarily elliptic in a compact set. It is however elliptic outside a large enough compact set since $ P_0 $ is close to $ P $,  or equivalently to  $ - \Delta $, at infinity and we may thus write
$$ P_0 = \widetilde{P}_0 + P_{\rm c} $$ with $ \widetilde{P}_0 $ uniformly elliptic and
$$ P_{\rm c} = \sum_{j,k=1}^d D_j \left( b_{jk}(x) D_k \right), \qquad b_{jk} \in C_0^{\infty}(\Ra^d) . $$
 We shall  absorb the contribution of $ \left( P_{\rm c} (hA+i)^{-1} u , (hA+i)^{-1} u \right) $ as in the original proof of Mourre \cite{Mourre},
by considering $ u$ which are spectrally localized very close to $ 0 $. Using (\ref{commutateurnonelliptique}) and the uniform ellipticity of $ \widetilde{P}_0 $ there exists $ c > 0 $ such that, for all $u$ satisfying (\ref{localisationspectrale}), we have
\begin{eqnarray}
 (i [P,F(hA)]u,u) & \geq & c h \big| \big| \nabla (hA+i)^{-1} u \big| \big|_{L^2}^2  - C h^2 \big| \big| |D| (hA+i)^{-1} u  \big| \big|_{L^2}^2 \nonumber \\ & & + \frac{h}{2} \left( P_{\rm c} (hA+i)^{-1} \varphi (P/\lambda) u , (hA+i)^{-1}  u \right) \label{calculexplicite} . 
\end{eqnarray}
Using (\ref{transformeedeRiesz}), we now introduce  
\begin{eqnarray*}
 R_{\rm c} =
 - \sum_{jk} R (j)^* b_{jk}(x) R(k) ,
\end{eqnarray*} 
ie $ R_{\rm c} = P^{-1/2} P_{\rm c} P^{-1/2} $ formally. Actually, by (\ref{proprietealgebriqueRiesz}), we have $ P^{1/2} R_{\rm c} P^{1/2} = P_{\rm c} $ at least in the form sense and this allows to rewrite  the last term of (\ref{calculexplicite}) as $ h/2 $ times the sum of the following two terms
\begin{eqnarray}
\left( R_{\rm c} \varphi (P/\lambda) P^{1/2} (hA+i)^{-1}  u , P^{1/2} (hA+i)^{-1} u \right) , \label{absorbeparcompacite} \\
\left( P_{\rm c}  [ (hA+i)^{-1} , \varphi (P/\lambda) ] u ,  (hA+i)^{-1}  u \right) .
\end{eqnarray}
The local compactness of the Riesz transform is crucial for the following result.
\begin{prop} \label{compaciteprop}  As $ \lambda \downarrow 0 $, we have 
$$ \big| \big| R_{\rm c} \varphi (P/\lambda) \big| \big|_{L^2 \rightarrow L^2} \rightarrow 0 . $$
\end{prop}

\noindent {\it Proof.} The operator $ R_{\rm c} \varphi (P/\lambda) $ can be written, for $ \lambda $ small enough, $ \left( R_{\rm c} \varphi (P) \right) \varphi (P/\lambda) $ since $ \varphi \equiv 1 $ near $ 0 $. The bracket is compact by Proposition \ref{conditionabstraite} and  $ \varphi (P/\lambda)  $ goes strongly to $ 0 $ as $ \lambda \downarrow 0 $, by the Spectral Theorem,
since $ 0 $ is not an eigenvalue of $ P $. Since $ R_{\rm c}  \varphi (P)$ is compact, $ \left( R_{\rm c} \varphi (P) \right) \varphi (P/\lambda) $ goes to $ 0 $ in operator norm. \finpreuve

\bigskip

By Proposition \ref{compaciteprop} and by choosing $ \lambda > 0 $ small enough, we can make  (\ref{absorbeparcompacite}) small so that, using (\ref{controlederivationbassefrequence}), we get the existence of $ c^{\prime} > 0 $ such that 
\begin{eqnarray}
 (i [P,F(hA)]u,u) & \geq & c^{\prime} h \big| \big| P^{1/2} (hA+i)^{-1} u \big| \big|_{L^2}^2  - C h^2 \big| \big| |D| (hA+i)^{-1} u  \big| \big|_{L^2}^2 \nonumber \\ & & - \frac{h}{2} \left| \big( P_{\rm c}[ (hA+i)^{-1}, \varphi (P/\lambda) ] u , (hA+i)^{-1}  u \big) \right| , 
 \label{compaciteexplicitee}
\end{eqnarray}
for all $ 0 < h < 1/4 $ and all $u$ satisfying (\ref{localisationspectrale}). It remains to deal with the last term of (\ref{compaciteexplicitee}). This is the purpose of the following proposition.



\begin{prop} For all $ \lambda > 0 $, there exists $ C_{\lambda} > 0 $ such that, for all $ v \in {\mathcal S}(\Ra^d) $ and all $ h$ 
$$  \big| \big| |D| [(hA+i)^{-1}, \varphi (P/\lambda)] v \big| \big|_{L^2 } \leq C_{\lambda} h \big| \big| P^{1/2} (hA+i)^{-1} v \big| \big|_{L^2}  . $$
\end{prop}

\noindent {\it Proof.} The proof relies on a standard combination of the resolvent identity
\begin{eqnarray}
 |D|[(hA+i)^{-1}, \varphi (P/\lambda)] = - h |D|(hA+i)^{-1} [ A, \varphi (P/\lambda)] (hA+i)^{-1} , \label{commutateurformal}
\end{eqnarray}
 and, for instance, the following Helffer-Sj\"ostrand formula (see \cite{DiSj1})
$$ \varphi (P/\lambda) =  \frac{1}{\pi} \int_{\Ca } \bar{\partial} \widetilde{\varphi}_{\lambda} (z) (P-z)^{-1} L(d z)  , $$
where $ L (dz) $ is the Lebesgue measure on $ \Ca \simeq \Ra^2 $ and $ \widetilde{\varphi}_{\lambda} \in C_0^{\infty} (\Ca) $ is an almost analytic extension of $ \varphi_{\lambda} := \varphi (\frac{\cdot}{\lambda}) $, ie which coincides with $ \varphi_{\lambda} $ on $ \Ra $ and such that $ \bar{\partial} \widetilde{\varphi}_{\lambda} = {\mathcal O} (|\mbox{Im}(z)|^{\infty}) $. We have
$$  [ A, \varphi (P/\lambda)] = - \frac{1}{\pi} \int_{\Ca} \bar{\partial} \widetilde{\varphi}_{\lambda} (z) (P-z)^{-1} [A,P] (P-z)^{-1} L(d z)  , $$
 hence, using (\ref{pourtraverserderivees}), we can rewrite (\ref{commutateurformal}) as
 $$  \frac{h}{\pi} \left( \int_{\Ca} \bar{\partial} \widetilde{\varphi}_{\lambda} (z) (hA+i+ih)^{-1} |D| P^{-1/2} \left\{ (P-z)^{-1}  P^{1/2} [A,P] P^{-1/2} (P-z)^{-1} \right\}  L(d z) \right) P^{1/2} (hA+i)^{-1} $$ 
where it is not hard to check that the operator $ \{ \ldots \} $ is bounded on $ L^2 $ with norm of polynomial growth with respect to $ |\mbox{Im}(z)|^{-1} $ (for $z$ in the support of $ \widetilde{\varphi}_{\lambda} $).
The result follows. \finpreuve

\bigskip

\noindent {\it End of the proof of Theorem \ref{theoremeconditionnel}.}
Since $ P_{\rm c} $ is of div-grad type, the last term of (\ref{compaciteexplicitee}) is bounded  by
$ - C_{\lambda} h^2 || P^{1/2} (hA+i)u||_{L^2}^2 $, from below. Thus, by choosing $ h $  and $ c_{\lambda} > 0 $ both small enough, we finally get
$$ (i [P,F(hA)]u,u)  \geq  c_{\lambda} h \big| \big| P^{1/2} (hA+i)^{-1} u \big| \big|_{L^2}^2  $$
for all $u$ satisfying (\ref{localisationspectrale}). We then obtain (\ref{estimeereduite}) as in the proof of Theorem  \ref{theoremeprincipal}. 
This completes the proof. \finpreuve


\end{document}